\documentclass[letterpaper]{amsproc}
\usepackage{xspace}

\newcommand{\baseRing}[1]{\ensuremath{\mathbb{#1}}}
\newcommand{\Z}{\baseRing{Z}}
\newcommand{\R}{\baseRing{R}}
\newcommand{\C}{\baseRing{C}}
\newcommand{\idM}{\baseRing{I}}
\newcommand{\pd}[2]{{\frac{\partial {#1}}{\partial {#2}}}}
\newcommand{\transpose}[1]{{#1}^t}
\newcommand{\CB}{\ensuremath{\mathcal{B}}\xspace}

\newcommand{\jgU}{\ensuremath{\mathfrak{U}}}

\theoremstyle{plain}
\newtheorem{theorem}{Theorem}[section]
\newtheorem{corollary}[theorem]{Corollary}
\newtheorem{prop}[theorem]{Proposition}

\theoremstyle{definition}
\newtheorem{definition}[theorem]{Definition}
\newtheorem{remark}[theorem]{Remark}
\newtheorem{example}[theorem]{Example}

\numberwithin{equation}{section}

\DeclareMathOperator{\aut}{Aut}

\DeclareMathOperator{\gr}{Gr}
\DeclareMathOperator{\vspan}{Span}
\DeclareMathOperator{\res}{Res}
\DeclareMathOperator{\homom}{Hom}
\DeclareMathOperator{\sym}{Sym}

\newcommand{\Script}[1]{\ensuremath{{\mathcal{#1}}}}
\newcommand{\HH}{\Script{H}}

\newcommand{\FF}{\Script{F}}

\newcommand{\KK}{\Script{K}}

\newcommand{\MM}{\Script{M}}
\newcommand{\QQ}{\Script{Q}}

\newcommand{\ScS}{\Script{S}}
\newcommand{\VV}{\Script{V}}
\newcommand{\scA}{\Script{A}}
\newcommand{\XX}{\Script{X}}
\newcommand{\CC}{\Script{C}}

\newcommand{\ie}{\textsl{i.e.}\xspace}
\newcommand{\del}{\ensuremath{\partial}}
\newcommand{\poly}{\ensuremath{(\Delta^*)^r}}

\newcommand{\ti}[1]{\tilde{#1}}
\newcommand{\conj}{\overline}

\newcommand{\GC}{\ensuremath{G_\C}}

\newcommand{\gl}{\ensuremath{\mathfrak{gl}}}
\newcommand{\jlg}{\ensuremath{\mathfrak{g}}}
\newcommand{\gp}{\ensuremath{\mathfrak{p}}}

\newcommand{\lgr}{\ensuremath{\mathfrak{g}_\R}}
\newcommand{\DD}{\ensuremath{D}}
\newcommand{\DC}{\ensuremath{\check{D}}}

\begin{document}

\title{Asymptotic Hodge Theory and Quantum Products}

\bibliographystyle{amsplain}

\address{Department of Mathematics and Statistics\\ University of
  Massachusetts\\ Amherst\\ MA 01003}

\author{Eduardo Cattani}
\author{Javier Fernandez}

\email{cattani@math.umass.edu}
\email{jfernand@math.umass.edu}

\begin{abstract}

Assuming suitable convergence properties for the Gromov-Witten
potential of a Calabi-Yau manifold $X$, one may construct a
polarized variation of Hodge structure over the complexified
K\"ahler cone of $X$.  
 In this paper we show that, in the case of 
  fourfolds, there is a correspondence between ``quantum potentials''
  and polarized variations of Hodge structures that degenerate to a
  maximally unipotent boundary point. Under this correspondence, the
WDVV equations are seen to be equivalent to the Griffiths' trasversality
property of a variation of Hodge structure.
\end{abstract}


\maketitle


\section{Introduction}
\label{sec:introduction}

The Gromov-Witten potential of a Calabi-Yau manifold is a generating
function for some of its enumerative data.  It may be written as
$\Phi^{GW} = \Phi_{0} + \Phi_{\hbar}$, where $\Phi_{0}$ is determined
by the cup product structure in cohomology.  For a quintic
hypersurface $X\subset P^4$, the quantum potential $\Phi_{\hbar}$ is
holomorphic in a neighborhood of $0\in \C$; the coefficients of its
expansion
\begin{equation*}
  \Phi_{\hbar}(q) \ :=\ \sum_{d=1}^\infty \langle
  I_{0,0,d}\rangle\cdot q^d,
\end{equation*}
are the Gromov-Witten invariants $\langle I_{0,0,d}\rangle$ which
encode information about the number of rational curves of degree $d$
in $X$.  This potential gives rise to a flat connection on the trivial
bundle over the punctured disk $\Delta^*$ with fiber $V = \oplus_p
H^{p,p}(X)$.  This flat bundle is shown to underlie a polarized
variation of Hodge structure whose degeneration at the origin is, in
an appropriate sense, maximal.  The Mirror Theorem in the context of
\cite{ar:CDGP-pair}, \cite{ar:LLY-mirror-1},
\cite{ar:givental-equivariant}, \cite[Theorem 11.1.1]{bo:CK-mirror}
asserts that this is the variation of Hodge structure arising from a
family of mirror Calabi-Yau threefolds and, therefore, that the
Gromov-Witten potential may be computed from the period map of this
family, written with respect to a canonical coordinate at
``infinity''.  This leads to the effective computation (and
prediction) of the number of rational curves, of a given degree, in a
quintic threefold.

The Mirror Theorem, in the sense sketched above, has a conjectural
generalization for toric Calabi-Yau threefolds \cite[Conjecture
8.6.10]{bo:CK-mirror} but the situation in the higher dimensional case
is considerably murkier.  Still, one may write a formal potential in
terms of axiomatically defined Gromov-Witten invariants and, assuming
suitable convergence properties, construct from it an abstract
polarized variation of Hodge structure.  The third derivatives of
$\Phi^{GW}$ may be used to define a quantum product on $H^*(X)$ whose
associativity is equivalent to a system of partial differential
equations satisfied by $\Phi^{GW}$.  These are the so-called WDVV
equations, after E.~Witten, R.~Dijkgraaf, H.~Verlinde, and
E.~Verlinde.

The purpose of this paper is to further explore the relationship
between quantum potentials and variations of Hodge structure.  The
weight-three case has been extensively considered in
\cite{bo:CK-mirror, ar:greg-higgs}; here we consider the case of
structures of weight four.

In \S \ref{sec:variations_at_infinity}, after recalling some basic
notions of Hodge theory, we describe the behavior of variations at
``infinity'' in terms of two ingredients: a nilpotent orbit and a
holomorphic function $\Gamma$ with values in a graded nilpotent Lie
algebra.  Nilpotent orbits may be characterized in terms of polarized
mixed Hodge structures; in the particular case when these split over
$\R$, their structure mimics that of the cohomology of a K\"ahler
manifold under multiplication by elements in the K\"ahler cone.  For a
Calabi-Yau fourfold $X$, this action ---together with the intersection
form--- characterizes the cup product in $\oplus_p H^{p,p}(X)$.  The
holomorphic function $\Gamma$, in turn, is completely determined by
one of its components, $\Gamma_{-1}$, relative to the Lie algebra
grading.  Moreover, it must satisfy the differential equation
(\ref{eq:integcond}) involving the monodromy of the variation.  This
is the content of Theorem~\ref{th:improved_2.8} which generalizes a
result of P.~Deligne \cite[Theorem 11]{ar:Del-local_behavior} for a
case when the variation of Hodge structure is also a variation of
mixed Hodge structure.  Throughout this section, and indeed this whole
paper, we restrict ourselves to real variations of Hodge structure
without any reference to integral structures.

The asymptotic data associated with a variation of Hodge structure
depends on the choice of local coordinates.  In
\S \ref{sec:canonical_coordinates} we describe how the nilpotent orbit
and $\Gamma$ behave under a change of coordinates and show that in
certain cases there are canonical choices of coordinates.  This
happens, for example, in the cases of interest in mirror symmetry and
we recover, in this manner, Deligne's Hodge theoretic description of
these canonical coordinates.

In \S \ref{sec:frobenius} we begin our discussion of quantum products
by concentrating on their ``constant'' part.  We define the notion of
polarized, graded Frobenius algebras and show that they give rise to
nilpotent orbits which are maximally degenerate in an appropriate
sense.  Moreover, in the weight-four case this correspondence may be
reversed.

Finally, in \S \ref{sec:quantum} we define a notion of quantum
potentials in polarized, graded Frobenius algebras of weight four
that abstracts the main properties of the Gromov-Witten potential for
Calabi-Yau fourfolds.  Following the arguments of \cite[\S
8.5.4]{bo:CK-mirror} we construct a variation of Hodge structure
associated with such a potential and explicitly describe its
asymptotic data.  We show, in particular, that the function $\Gamma$
is completely determined by the potential in a manner that transforms
the WDVV equations into the horizontality equation
(\ref{eq:integcond}).  This is done in Theorem~\ref{main} which
establishes an equivalence between such potentials and variations with
certain prescribed limiting behavior.  For the case of weight three a
similar theorem is due to G.~Pearlstein \cite[Theorem
7.20]{ar:greg-higgs}.

A distinguishing property of the weight $3$ and $4$ cases is that the
quantum product on the $(p,p)$-cohomology is determined by the
$H^{1,1}$-module structure. In the case of constant products, this
statement is the content of
Proposition~\ref{prop:constant_prod_from_nilpotent}.  For general
weights the PVHS on the complexified K\"ahler cone determines only the
structure of the $(p,p)$-cohomology as a $\sym H^{1,1}$-module,
relative to the small quantum product. At the abstract level, one can
prove that there is a correspondence between germs of PVHS of weight
$k$ at a maximally unipotent boundary point whose limiting mixed Hodge
structure is Hodge-Tate and flat, real families of polarized graded
Frobenius $\sym H^{1,1}$-modules over $\oplus_{p=0}^k H^{p,p}$. 
Details will appear elsewhere.

Barannikov (\cite{ar:barannikov-1},~\cite{ar:barannikov-2} and
\cite{ar:barannikov-3}) has shown that for projective complete
intersections, the PVHS constructed from the Gromov-Witten potential
is of geometric origin and coincides with the variation arising from
the mirror family. He has introduced, moreover, the notion of
semi-infinite variations of Hodge structure. These more general 
variations are shown to correspond to solutions of the WDVV equations.

\smallskip

\noindent{\bf Acknowledgments:}  We are very grateful to Emma Previato
for organizing the Special Session on Enumerative Geometry in Physics
at the Regional AMS meeting of April 2000 in Lowell, Massachusetts and
for putting together this volume.
We also wish to thank David Cox for very helpful conversations and to 
Serguei Barannikov for his useful comments.


\section{Variations at Infinity}
\label{sec:variations_at_infinity}

We begin by reviewing some basic results about the asymptotic behavior
of variations of Hodge structure.  We refer to \cite{ar:CK-luminy,
  ar:Gri-periods-1, bo:griffiths-topics, ar:S-vhs} for details and
proofs.

Let $M$ be a connected complex manifold, a (real) {\em variation of
  Hodge structure} (VHS) $\VV$ over $M$ consists of a holomorphic
vector bundle $\VV\to M$, endowed with a flat connection $\nabla$, a
flat real form $\VV_{\R}$, and a finite decreasing filtration $\FF$ of
$\VV$ by holomorphic subbundles ---the \textsl{ Hodge filtration\/}---
satisfying
\begin{eqnarray}\label{horizontality}
  \nabla\FF^p & \subset & \Omega^1_M
  \otimes \FF^{p-1}\quad\hbox{(Griffiths' horizontality) and}
  \\ \label{opposed}
  \VV& = & \FF^p \oplus \conj{\FF}^{k-p+1}
\end{eqnarray}
for some integer $k$ ---the \textsl{weight\/} of the variation--- and
where $\conj{\FF}$ denotes conjugation relative to $\VV_{\R}$.  As a
$C^{\infty}$-bundle, $\VV$ may then be written as a direct sum
\begin{equation}\label{bundledecomposition}
  \VV = \bigoplus_{p+q=k}\ \VV^{p,q}\ ,\quad \quad \VV^{p,q} = \FF^p
  \cap \conj{\FF}^q\,;
\end{equation}
the integers $h^{p,q} = \dim\,\VV^{p,q}$ are the \textsl{ Hodge
  numbers\/}.  A \textsl{ polarization\/} of the VHS is a flat
non-degenerate bilinear form $\ScS$ on $\VV$, defined over $\R$, of
parity $(-1)^k$, whose associated flat Hermitian form
$\ScS^h(\,.\,,\,.\,) = i^{-k}\, \ScS(\,.\,,\,\bar.\,)$ is such that
the decomposition~(\ref{bundledecomposition}) is $\ScS^h$-orthogonal
and $(-1)^p\ScS^h$ is positive definite on $\VV^{p,k-p}$.

Specialization to a fiber defines the notion of \textsl{polarized
  Hodge structure} on a $\C$-vector space $V$.  We will denote by $D$
\cite{ar:Gri-periods-1} the \textsl{ classifying space\/} of all
polarized Hodge structures of given weight and Hodge numbers on a
fixed vector space $V$, endowed with a fixed real structure $V_{\R}$
and the polarizing form $S$.  Its Zariski closure $\check D$ in the
appropriate variety of flags consists of all filtrations $F$ in $V$,
with $\dim\,F^p = \sum_{r\geq p}\,h^{r,k-r}$, satisfying
$\,S(F^p,F^{k-p+1}) = 0\, $.  The complex Lie group $G_{\C}=
\aut(V,S)$ acts transitively on $\check D$ ---therefore $\check D$ is
smooth--- and the group of real points $G_{\R}$ has $D$ as an open
dense orbit.  We denote by $\jlg \subset \gl(V)$, the Lie algebra of
$G_{\C}$, and by $\lgr \subset \jlg$ that of $G_{\R}$.  The choice
of a base point $F\in \check D$ defines a filtration in $\jlg$
\begin{equation*}
  F^a\jlg = \{\,T\in\jlg\ :\ T\,F^p \subset F^{p+a}\,\}\,.
\end{equation*}
The Lie algebra of the isotropy subgroup $B\subset G_{\C}$ at $F$ is
$F^0\jlg$ and $F^{-1}\jlg/F^0\jlg$ is an ${\rm Ad}(B)$-invariant
subspace of $\jlg /F^0\jlg$.  The corresponding $G_{\C}$-invariant
subbundle of the holomorphic tangent bundle of $\check D$ is the
\textsl{ horizontal tangent bundle\/}.  A polarized VHS over a
manifold $M$ determines ---via parallel translation to a typical fiber
$V$--- a holomorphic map $\Phi\colon M \to D/\Gamma$ where $\Gamma$ is
the monodromy group (Griffiths' period map).  By definition, it has
\textsl{horizontal} local liftings into $D$, \ie, its differentials
take values on the horizontal tangent bundle.

\begin{example}
  Let $X$ be an $n$-dimensional, smooth projective variety, $\omega\in
  H^{1,1}(X)$ a K\"ahler class.  For any $k=0,\dots,2n$, the Hodge
  decomposition (see~\cite{bo:griffiths-principlesAG})
  \begin{equation*}
    H^k(X,\C) \ = \ \bigoplus_{p+q =k} H^{p,q}(X); \quad \conj{H^{p,q}}
    = H^{q,p}
  \end{equation*}
  determines a Hodge structure of weight $k$ by
  \begin{equation*}
    F^{p} \ := \ \bigoplus_{a\geq p} H^{a,k-a}.
  \end{equation*}
  Its restriction to the \textsl{primitive cohomology}
  \begin{equation*}
    H_0^{n-\ell}(X,\C) \ :=\ \{\alpha \in H^{n-\ell}(X,\C)\,:\,
    \omega^{\ell+1} \cup \alpha = 0\}\,,\ \ell\geq 0,
  \end{equation*}
  is polarized by the form $Q_\ell(\alpha ,\beta) = Q(\alpha ,
  \beta\cup \omega^\ell)$, $\alpha, \beta\in H_0^{n-\ell}(X,\C)$, and
  where $Q$ denotes the signed intersection form given, for $\alpha\in
  H^k(X,\C)$, $\alpha'\in H^{k'}(X,\C)$ by:
  \begin{equation}\label{intersectionform}
    Q(\alpha,\alpha') \ :=\ (-1)^{k(k-1)/2}\ \int_X \alpha\cup\alpha'.
  \end{equation}
  A family $\XX \to M$ of smooth projective varieties gives rise to a
  polarized variation of Hodge structure $\VV \to M$, where $\VV_m
  \cong H_0^k(\XX_m,\C)$, $m\in M$.
\end{example}

Our main concern is the asymptotic behavior of $\Phi$ near the
boundary of $M$, with respect to some compactification $\conj{M}$
where $\conj{M}-M$ is a divisor with normal crossings (the divisor at
``infinity").  Such compactifications exist, for instance, if $M$ is
quasiprojective.  Near a boundary point $p\in\conj{M}-M$ we can choose
an open set $W$ such that $W\cap M \simeq (\Delta^*)^r \times
\Delta^m$ and then consider the local period map
\begin{equation}\label{localperiod}
  \Phi:(\Delta^*)^r \times \Delta^m\rightarrow \DD/\Gamma.
\end{equation}
We shall also denote by $\Phi$ its lifting to the universal covering
$U^r \times \Delta^m$, where $U$ denotes the upper-half plane.  We
denote by $z=(z_j)$, $t=(t_l)$ and $s=(s_j)$ the coordinates on $U^r$,
$\Delta^m$ and $(\Delta^*)^r$ respectively. By definition, we have
$s_j = e^{2\pi i z_j}$.

According to Schmid's Nilpotent Orbit Theorem \cite{ar:S-vhs}, the
singularities of $\Phi$ at the origin are, at worst, logarithmic; this
is essentially equivalent to the regularity of the connection
$\nabla$.  More precisely, assuming quasi-unipotency ---this is
automatic in the geometric case--- and after passing, if necessary, to
a finite cover of $(\Delta^*)^r$, there exist commuting nilpotent
elements\footnote{Our sign convention is consistent with
  \cite{ar:Del-local_behavior, ar:S-vhs} but opposite to that in
  \cite{bo:CK-mirror, ar:morrison-guide}.} $N_1,\dots,N_r\in \jlg_\R$,
with $N^{k+1}=0$ and such that
\begin{equation}\label{not}
  \Phi(s,t) \ =\ \exp\left(\sum_{j=1}^r\frac {\log s_j}{2\pi i}
    N_j\right)\cdot \Psi(s,t)\,,
\end{equation}
where $\Psi \colon \Delta^{r+m} \to \check D$ is holomorphic.

We will refer to $N_1,\dots,N_r$ as the \textsl{ local monodromy
  logarithms} and to $F_0:=\Psi(0)$ as the \textsl{limiting Hodge
  filtration}.  They combine to define a \textsl{nilpotent orbit}
$\{N_1,\dots,N_r; F_0\}$.  We recall:

\begin{definition} With notation as above,
  $\{N_1,\dots,N_r; F_0\}$ is called a nilpotent orbit if the map
  \begin{equation*}
    \theta(z) = \exp(\sum_{j=1}^r z_j
    N_j)\cdot F_0
  \end{equation*}
  is horizontal and there exists $\alpha\in \R$ such that $\theta(z)
  \in D$ for ${\rm Im}(z_j) > \alpha$.
\end{definition}
Theorem~\ref{th:nilporbit} below gives an algebraic characterization
of nilpotent orbits which will play a central role in the sequel.

We point out that the local monodromy is topological in nature, while
the limiting Hodge filtration depends on the choice of coordinates
$s_j$.  To see this, we consider, for simplicity, the case $m=0$.  A
change of coordinates compatible with the divisor structure must be,
after relabeling if necessary, of the form $(s'_1,\dots,s'_r) =
(s_1f_1(s),\dots,s_rf_r(s))$ where $f_j$ are holomorphic around $0\in
\Delta^r$ and $f_j(0) \not=0$. We then have from~(\ref{not}),
\begin{equation}\label{psicoord}
  \begin{split}
    \Psi'(s') &= \exp(- \frac{1}{2\pi i} \sum_{j=1}^r
    \log(s'_j)N_j)\cdot \Phi(s')\\ &= \exp(- \frac{1}{2\pi i}
    \sum_{j=1}^r \log(f_j)N_j) \exp(- \frac{1}{2\pi i} \sum_{j=1}^r
    \log(s_j)N_j) \cdot \Phi(s)\\ &= \exp(- \frac{1}{2\pi i}
    \sum_{j=1}^r \log(f_j)N_j) \cdot \Psi(s),
  \end{split}
\end{equation}
and, letting $s\to 0$
\begin{equation}
  \label{eq:change_F}
  F'_0 = \exp (- \frac {1}{2\pi i} \sum_j \log (f_j(0)) N_j) \cdot
  F_0.
\end{equation}

These constructions may also be understood in terms of Deligne's
canonical extension \cite{bo:Del-equations}.  Let $\VV \to
(\Delta^*)^r \times \Delta^m$ be the local system underlying a
polarized VHS and pick a base point $(s_0,t_0)$.  Given $v\in V :=
\VV_{(s_0,t_0)}$, let $v^{\flat}$ denote the multivalued flat section
of $\VV$ defined by $v$.  Then
\begin{equation}\label{cansections}
  \tilde v(s,t) \ :=\ \exp\left(\sum\frac {\log s_j}{2\pi i}
    N_j\right)\cdot v^{\flat}(s,t)
\end{equation}
is a global section of $\VV$. The canonical extension $\conj\VV \to
\Delta^{r+m}$ is characterized by its being trivialized by sections of
the form (\ref{cansections}).  The Nilpotent Orbit Theorem then
implies that the Hodge bundles $\FF^p$ extend to holomorphic
subbundles $\conj\FF^p\subset \conj\VV$.  Writing the Hodge bundles in
terms of a basis of sections of the form (\ref{cansections}) yields
the holomorphic map $\Psi$.  Its constant part ---corresponding to the
nilpotent orbit--- defines a polarized VHS as well.

A nilpotent linear transformation $N\in\gl(V_\R)$ defines an
increasing filtration, the \textsl{weight filtration}, $W(N) $ of $V$,
defined over $\R$ and uniquely characterized by requiring that
$N(W_l(N))\subset W_{l-2}(N)$ and that $N^l:\gr_{l}^{W(N)}\rightarrow
\gr_{-l}^{W(N)}$ be an isomorphism.  It follows from
\cite[Theorem~3.3]{ar:CK-polarized} that if $N_1,\dots,N_r$ are local
monodromy logarithms arising from a polarized VHS then the weight
filtration $W(\sum \lambda_j N_j)$, $\lambda_j\in \R_{>0}$, is
independent of the choice of $\lambda_1,\dots,\lambda_r$ and,
therefore, is associated with the positive real cone $\CC \subset
\jlg_\R$ spanned by $N_1,\dots,N_r$.

A \textsl{mixed Hodge Structure} (MHS) on $V$ consists of a pair of
filtrations of $V$, $(W, F)$, $W$ defined over $\R$ and increasing,
$F$ decreasing, such that $F$ induces a Hodge structure of weight $k$
on $\gr_k^{W}$ for each $k$.  Equivalently, a MHS on $V$ is a
bigrading
\begin{equation*}
  V\ = \ \bigoplus I^{p,q}
\end{equation*}
satisfying $I^{p,q}\equiv \conj{I^{q,p}} \mod(\oplus_{a<p,b<q}
I^{a,b})$ (see~\cite[Theorem~2.13]{ar:CKS}).  Given such a bigrading
we define: $ W_l = \oplus_{p+q \leq l} I^{p,q}$, $F^a = \oplus_{p\geq
  a} I^{p,q}$.  A MHS is said to \textsl{split} over $\R$ if $I^{p,q}=
\conj{I^{q,p}}$; in that case the subspaces $V_l = \oplus_{p+q = l}
I^{p,q}$ define a real grading of $W$.  A map $T \in \gl(V)$ such that
$T(I^{p,q}) \subset I^{p+a,q+b}$ is called a morphism of bidegree
$(a,b)$.

A \textsl{polarized MHS} (PMHS)~\cite[(2.4)]{ar:CK-polarized} of
weight $k$ on $V_\R$ consists of a MHS $(W,F)$ on $V$, a $(-1,-1)$
morphism $N\in \lgr$, and a nondegenerate bilinear form $Q$ such that
\begin{enumerate}
\item $N^{k+1}=0$,
\item $W = W(N)[-k]$, where $W[-k]_j = W_{j-k}$,
\item $Q(F^a,F^{k-a+1}) = 0$ and,
\item the Hodge structure of weight $k+l$ induced by $F$ on
  $\ker(N^{l+1}:\gr_{k+l}^{W}\rightarrow \gr_{k-l-2}^{W})$ is
  polarized by $Q(\cdot,N^l \cdot)$.
\end{enumerate}

\begin{theorem} \label{th:nilporbit}
  Given a nilpotent orbit $\theta(z) = \exp(\sum_{j=1}^r z_j N_j)\cdot
  F $, the pair $(W(\CC),F)$ defines a MHS polarized by every $N\in
  \CC$.  Conversely, given commuting nilpotent elements
  $\{N_1,\ldots,N_r\}\in \lgr$ with the property that the weight
  filtration $W(\sum \lambda_j N_j)$, $\lambda_j\in \R_{>0}$, is
  independent of the choice of $\lambda_1,\dots,\lambda_r$, \ if $F
  \in \DC$ is such that $(W(\CC),F)$ is polarized by every element
  $N\in \CC$, then the map $\theta(z) = \exp(\sum_{j=1}^r z_j
  N_j)\cdot F$ is a nilpotent orbit.  Moreover, if $(W(\CC),F)$ splits
  over $\R$, then $\theta(z) \in D$ for ${\rm Im}(z_j) > 0$.
\end{theorem}

The first part of Theorem~\ref{th:nilporbit} was proved by Schmid
\cite[Theorem~6.16]{ar:S-vhs} as a consequence of his $SL_2$-orbit
theorem.  The converse is Proposition~4.66 in \cite{ar:CKS}.  The
final assertion is a consequence of
\cite[Proposition~2.18]{ar:CK-polarized}.

\begin{example}\label{ex:totalcohomology}
  Let $X$ be an $n$-dimensional, smooth projective variety.  Let $V =
  H^*(X,\C)$, $V_\R = H^*(X,\R)$.  The bigrading $I^{p,q} :=
  H^{n-p,n-q}(X)$ defines a MHS on $V$ which splits over $\R$.  The
  weight and Hodge filtrations are then
  \begin{equation*}
    W_l \ =\ \bigoplus_{d\geq 2n-l} H^d(X,\C),\quad F^p \ =\
    \bigoplus_{s}\bigoplus_{r\leq n-p}H^{r,s}(X).
  \end{equation*}
  Given a K\"ahler class $\omega\in H^{1,1}(X,\R):= H^{1,1}(X) \cap
  H^2(X,\R)$, let $L_{\omega}\in \gl(V_\R)$ denote multiplication by
  $\omega$.  Note that $L_{\omega}$ is an infinitesimal automorphism
  of the form (\ref{intersectionform}) and is a $(-1,-1)$ morphism of
  $(W,F)$.  Moreover, the Hard Lefschetz Theorem and the Riemann
  bilinear relations are equivalent to the assertion that $L_\omega$
  polarizes $(W,F)$.  Let $\KK \subset H^{1,1}(X,\R)$ denote the
  K\"ahler cone and
  \begin{equation*}
    \KK_\C := H^{1,1}(X,\R) \oplus i \KK \subset H^2(X,\C)
  \end{equation*}
  the complexified K\"ahler cone.  It then follows from
  Theorem~\ref{th:nilporbit} that for every $\xi\in \KK_\C$, the
  filtration $\exp(L_\xi)\cdot F$ is a Hodge structure of weight $n$
  on $V$ polarized by $Q$.  The map $\, \xi\in \KK_\C \mapsto
  \exp(L_\xi)\cdot F \,$ is the period map (in fact, the nilpotent
  orbit) of a variation of Hodge structure over $\KK_\C$.  Note that
  we can restrict the above construction to $V = \oplus_p H^{p,p}(X)$;
  this is the case of interest in mirror symmetry.
\end{example}

\begin{remark}
  The notion of nilpotent orbit is closely related to that of
  Lefschetz modules introduced by Looijenga and
  Lunts~\cite{ar:LL-lefschetz_modules}.  Indeed, it follows from
  \cite[Proposition~1.6]{ar:LL-lefschetz_modules} that if $(W,F)$ is a
  MHS, polarized by every $N$ in a cone $\CC$, and $\scA$ denotes the
  linear span of $\CC$ in $\lgr$, then $V$ is a Lefschetz module of
  $\scA$.
\end{remark}

Let now $\Phi$ be as in (\ref{localperiod}), $\{N_1,\dots,N_r;F_0\}$
the associated nilpotent orbit and $(W(\CC),F_0)$ the \textsl{limiting
  mixed Hodge structure}.  The bigrading $I^{*,*}$ of $V$ defined by
$(W(\CC),F_0)$ defines a bigrading $I^{*,*}\jlg$ of the Lie algebra
$\jlg$ associated with the MHS $(W(\CC)\jlg,F_0\jlg)$.  Set
\begin{equation}
  \label{eq:def_pa}
  \gp_a \ := \ \bigoplus_{q}I^{a,q}\jlg \quad \text{ and }\quad \jlg_-
  \ := \ \bigoplus_{a\leq
    -1}\gp_a.
\end{equation}
The nilpotent subalgebra $\jlg_-$ is a complement of the stabilizer
subalgebra at $F_0$. Hence $(\jlg_-, X \mapsto \exp(X)\cdot F_0)$
provides a local model for the $\GC$-homogeneous space $\DC$ near
$F_0$ and we can rewrite (\ref{not}) as:
\begin{equation}\label{gamma}
  \Phi(s,t) \ =\ \exp\left(\sum_{j=1}^r\frac {\log s_j}{2\pi i}
    N_j\right)\ \exp \Gamma(s,t)\cdot F_0\,,
\end{equation}
where $\Gamma\colon \Delta^{r+m} \to \jlg_-$ is holomorphic and
$\Gamma(0) = 0$.  The lifting of $\Phi$ to $U^r\times \Delta^m$ may
then be expressed as:
\begin{equation*}
  \Phi(z,t) = \exp X(z,t)\cdot F_0
\end{equation*}
with $X:U^r\times \Delta^m\rightarrow \jlg_-$ holomorphic.  Setting
$E(z,t) := \exp X(z,t)$, the horizontality of $\Phi$ is then expressed
by:
\begin{equation}
  \label{eq:horiz}
  E^{-1}\, dE = dX_{-1} \in \gp_{-1}\otimes T^*(U^r\times \Delta^m).
\end{equation}
The following explicit description of period mappings near infinity is
proved in ~\cite[Theorem~2.8]{ar:CK-luminy}.

\begin{theorem}\label{th:2.8}
  Let $\{N_1,\ldots,N_r;F_0\}$ be a nilpotent orbit and
  $\Gamma:\Delta^r\times \Delta^m \rightarrow \jlg_-$ be holomorphic,
  such that $\Gamma(0,0) = 0$. If the map
  \begin{equation}
    \label{eq:2.8_form}
    \Phi(z,t)=\exp(\sum z_j N_j)\ \exp(\Gamma(s,t))\cdot F_0
  \end{equation}
  is horizontal (\ie, (\ref{eq:horiz}) is satisfied), then $\Phi(z,t)$
  is a period mapping.
\end{theorem}

Given a period mapping as in (\ref{eq:2.8_form}), let $\Gamma_{-1}(s,t)$
denote the $\gp_{-1}$-component of $\Gamma$.  Then,
\begin{equation*}
  X_{-1}(z,t)\ =\ \sum_{j=1}^r z_j N_j+\Gamma_{-1}(s,t)\,, \text{ with }
  s_j=e^{2\pi i z_j},
\end{equation*}
and it follows from (\ref{eq:horiz}) that
\begin{equation} \label{eq:integcond}
  dX_{-1}\wedge dX_{-1}=0
\end{equation}
The following theorem shows that this equation characterizes period
mappings with a given nilpotent orbit.

\begin{theorem}
  \label{th:improved_2.8}
  Let $R:\Delta^r \times \Delta^m \rightarrow \gp_{-1}$ be a
  holomorphic map with $R(0)=0$.  Let $X_{-1}(z,t) = \sum_{j=1}^r z_j
  N_j+R(s,t)$, $s_j = e^{2\pi i z_j}$, and suppose that the
  differential equation (\ref{eq:integcond}) holds.  Then, there
  exists a unique period mapping (\ref{gamma}) defined in a
  neighborhood of the origin in $\Delta^{r+m}$ and such that
  $\Gamma_{-1} = R$.
\end{theorem}

\begin{proof}
  To prove uniqueness, we begin by observing that if $\Phi$ and
  $\Phi'$ are period mappings with the same associated nilpotent orbit
  and $\Gamma_{-1}(s,t) = \Gamma_{-1}'(s,t)$ then, for any $v \in
  F^p$, we may consider the sections $\nu(s,t) = E(s,t)\cdot
  v^\flat(s,t)$ and $\nu'(s,t) = E'(s,t)\cdot v^\flat(s,t)$ of the
  canonical extension $\conj{\mathcal V}$.  Clearly, $\nu(s,t)\in
  {\mathcal F}_{(s,t)}^p$ and $\nu'(s,t)\in {\mathcal F'}_{(s,t)}^p$.
  On the other hand, since $\Gamma_{-1}(s,t) = \Gamma_{-1}'(s,t)$, it
  follows that $E_{-1}(s,t) = E_{-1}'(s,t)$ and, consequently,
  $\nu(s,t)-\nu'(s,t)$ is a $\nabla$-flat section which extends to the
  origin and takes the value zero there.  Hence, $\nu(s,t)-\nu'(s,t)$
  is identically zero and ${\mathcal F}_{(s,t)}^p = ({\mathcal
    F}_{(s,t)}')^p$ for all values of $(s,t)$.

  To complete the proof of the Theorem it remains to show the
  existence of a period mapping with given nilpotent orbit and
  $\Gamma_{-1}(s,t) = R(s,t)$.  This amounts to finding a solution to
  the differential equation (\ref{eq:horiz}) with
  \begin{equation*}
    X_{-1}(z,t)\ =\ \sum_{j=1}^r z_j N_j+ R(s,t)\,, \text{ with }
    s_j=e^{2\pi i z_j},
  \end{equation*}
  assuming that the integrability condition (\ref{eq:integcond}) is
  satisfied.  Set $G(s,t) = \exp \Gamma(s,t)$ and $\Theta= d
  (\sum_{j=1}^r z_j N_j)$.  Then (\ref{eq:horiz}) may be rewritten as
  \begin{equation}
    \label{eq:e5}
    d G = [G,\Theta] + G d G_{-1} \text{ with } G(0,0)=\idM,
  \end{equation}
  where $\idM$ denotes the identity, while the
  condition~(\ref{eq:integcond}) takes the form:
  \begin{equation}
    \label{eq:e6}
    d G_{-1} \wedge \Theta + \Theta \wedge d G_{-1} + d G_{-1}\wedge d
    G_{-1} = 0.
  \end{equation}

  By considering the $\gp_{-l}$-graded components of~(\ref{eq:e5}) we
  obtain a sequence of equations:
  \begin{equation}
    \label{eq:e7}
    d G_{-l} = [G_{-l+1}, \Theta] + G_{-l+1} d G_{-1}, \quad
    G_{-l}(0,0)=0, \quad l\geq 2.
  \end{equation}

  Assume inductively that, for $l\geq 2$, we have constructed
  $G_{-l+1}$ satisfying~(\ref{eq:e7}) and such that
  \begin{equation*}
    d G_{-l+1} \wedge \Theta + \Theta \wedge d G_{-l+1} + d
    G_{-l+1}\wedge d G_{-1} = 0.
  \end{equation*}
  Then, the initial value problem
  \begin{equation*}
    d G_{-l} = [G_{-l+1}, \Theta] + G_{-l+1} d G_{-1},\quad
    G_{-l}(0,0) = 0,
  \end{equation*}
  has a solution which verifies
  \begin{eqnarray*}
    \lefteqn{d G_{-l} \wedge \Theta + \Theta \wedge d G_{-l} + d
      G_{-l} \wedge d G_{-1} =}\\ &=& [ G_{-l+1}, \Theta] \wedge
    \Theta + G_{-l+1} d G_{-1} \wedge \Theta + \Theta \wedge
    [G_{-l+1}, \Theta] + \Theta \wedge G_{-l+1} d G_{-1} + \\ & &
    \mbox{} + [G_{-l+1}, \Theta] \wedge d G_{-1} + G_{-l+1} d G_{-1}
    \wedge d G_{-1} \\ &=& -\Theta \wedge G_{-l+1} \Theta +G_{-l+1} d
    G_{-1} \wedge \Theta + \Theta \wedge G_{-l+1} \Theta + \Theta
    \wedge G_{-l+1} d G_{-1} +\\ & & \mbox{} + G_{-l+1} \Theta \wedge
    d G_{-1} + G_{-l+1} d G_{-1} \wedge d G_{-1} - \Theta \wedge
    G_{-l+1} d G_{-1} \\ &=& G_{-l+1} (d G_{-1} \wedge \Theta + \Theta
    \wedge d G_{-1} + d G_{-1} \wedge d G_{-1}) = 0.
  \end{eqnarray*}
  Thus we may, inductively, construct a solution of
  ~(\ref{eq:e6}). Theorem~\ref{th:2.8} now implies that the map
  \begin{equation*}
    \Phi(z,t) = \exp(\sum_{j=1}^r z_j N_j)\ G(s,t)\cdot F_0
  \end{equation*}
  is the desired period map.
\end{proof}

\begin{remark}
  The uniqueness part of the argument is contained in Lemmas~2.8
  and~2.9 of \cite{ar:CDK}, while the existence proof is contained in
  the unpublished manuscript \cite{prep:C-addenda}.  A particular case
  of Theorem~\ref{th:improved_2.8} is given in
  \cite[Theorem~11]{ar:Del-local_behavior}; a generalization to the
  case of variations of MHS appears in \cite{ar:greg-higgs}.
\end{remark}


\section{Canonical Coordinates}
\label{sec:canonical_coordinates}

The asymptotic data of a polarized variation of Hodge structure over
an open set $W \simeq (\Delta^*)^r \times \Delta^m$ depends on the
choice of coordinates on the base.  We have already observed that the
local monodromy logarithms $N_j$ are independent of coordinates and
have shown in (\ref{eq:change_F}) how the limiting Hodge filtration
changes under a coordinate transformation. Here we will discuss the
dependence of the holomorphic function $\Gamma \colon (\Delta^*)^r
\times \Delta^m \to \gp_{-1}$ and show that, in special cases, there
is a natural choice of coordinates.  This choice will be seen to agree
with that appearing in the mirror symmetry setup and which has already
been given a Hodge-theoretic interpretation by Deligne
\cite{ar:Del-local_behavior}.  These canonical coordinates may also be
interpreted, in the case of families of Calabi-Yau threefolds as the
coordinates where the Picard-Fuchs equations take on a certain
particularly simple form (\cite[Prop. 5.6.1]{bo:CK-mirror}).  To
simplify the discussion we will restrict our discussion to the case
$m=0$.  The general case, which follows easily, will not be needed in
the sequel.

Since we are required to preserve the divisor structure at the boundary,
we want to study the behavior of the asymptotic data under coordinate
changes of the form
\begin{equation}
  \label{eq:coordinate_transformation}
  s'_j = s_j f_j(s)
\end{equation}
where the functions $f_j$ are holomorphic  in a neighborhood of $0
\in \Delta^r$ and $f_j(0)\neq 0$.

Given a PVHS over $(\Delta^*)^r$ and a choice of local coordinates
$(s_1,\dots,s_r)$  around
$0$, we write the associated period map as in (\ref{gamma}):
\begin{equation*}
  \Phi(s) \ =\
  \exp(\sum_{j=1}^r \frac{\log s_j}{2\pi i} N_j)
  \exp(\Gamma(s)) \cdot F_0
\end{equation*}
Given another system of coordinates $s' = (s'_1,\dots,s'_r)$ as in
(\ref{eq:coordinate_transformation}), let $F'_0$ and $\Gamma'$ denote
the corresponding asymptotic data.  By (\ref{eq:change_F}), $F'_0 =
{\mathcal M}\cdot F_0$, where
\begin{equation*}
   {\mathcal M}\ := \ \exp(-\frac{1}{2\pi
   i}\sum_{j=1}^r\log f_j(0)\,N_j).
\end{equation*}

\begin{prop}\label{prop:change_coordinates}
  Under a coordinate change as in
  (\ref{eq:coordinate_transformation}):
  \begin{equation}
    \label{eq:change_gamma}
    {\mathcal M}^{-1}\exp(\Gamma'(s')) {\mathcal M} \ =\
    \exp\left(-\frac{1}{2\pi i} \sum_{j=1}^r \log
        \frac{f_j(s)}{f_j(0)}\,N_j\right) \exp(\Gamma(s)).
  \end{equation}
\end{prop}

\begin{proof}  Let $W$ denote the filtration $W(N)[-k]$, where
  $N$ is an arbitrary element in the cone $\CC$ positively spanned by
  $N_1,\dots,N_r$.  Note that $\MM$ leaves $W$ invariant.  Moreover,
  since the monodromy logarithms are $(-1,-1)$-morphisms of the mixed
  Hodge structure $(W,F_0)$ it follows easily that
  \begin{equation*}
    I^{a,b}(W,F'_0) \ =\ \MM\cdot I^{a,b}(W,F_0),
  \end{equation*}
  where $I^{*,*}(W,F_0)$ denotes the canonical bigrading of the MHS.
  This implies that the associated bigrading of the Lie algebra $\jlg$
  and, in particular, that the subalgebra $\jlg_-$ defined in
  (\ref{eq:def_pa}) are independent of the choice of coordinates.

  According to (\ref{psicoord}), $\Psi'(s') \ =\ \exp\left(\frac
    1{2\pi i}\sum_j {\log f_j(s)} N_j\right)\cdot \Psi(s)$, therefore
  \begin{equation}\label{intermediate}
    \exp \Gamma'(s')\, \MM\cdot F_0 \ =\ \exp\left(\frac {1}{2\pi
    i}\sum_{j=1}^r \log f_j(s)\, N_j\right)\ \exp \Gamma(s)\cdot F_0.
  \end{equation}
  This identity, in turn, implies (\ref{eq:change_gamma}) since the
  group elements in both sides of (\ref{intermediate}) lie in
  $\exp(\jlg_-)$.
\end{proof}

\begin{corollary}\label{cor:change_coordinates-1}
  With the same notation of Proposition~\ref{prop:change_coordinates},
  \begin{equation}
    \label{eq:change_w1}
    {\mathcal M}^{-1}\Gamma'_{-1}{\mathcal M}\ = \ -\frac{1}{2\pi i}
    \sum_j \log  \frac{f_j(s)}{f_j(0)}\, N_j + \Gamma_{-1}.
  \end{equation}
\end{corollary}

\begin{proof}
  This follows considering the $\gp_{-1}$-component
  in~(\ref{eq:change_gamma}), given the observation that this subspace
  is invariant under coordinate changes.
\end{proof}

Up to rescaling, we may assume that our coordinate change
(\ref{eq:coordinate_transformation}) satisfies $f_j(0) =1$,
$j=1\dots,r$.  Such changes will be called {\sl simple}.  In this case
$\MM = \idM$, $F'_0 = F_0$, and the transformation
(\ref{eq:change_w1}) is just a translation in the direction of the
nilpotent elements $N_j$.  Thus, whenever the subspace spanned by
$N_1,\dots,N_r$ has a natural complement in $\gp_{-1}$ we will be able
to choose coordinates, unique up to scaling, such that $\Gamma_{-1}$
takes values in that complement.  This is the situation in the
variations of Hodge structure studied in mirror symmetry.  In this
context, one analyzes the behavior of PVHS near some special boundary
points.  They come under the name of ``large radius limit points''
(see~\cite[\S 6.2.1]{bo:CK-mirror}) or ``maximally unipotent boundary
points'' (see~\cite[\S 5.2]{bo:CK-mirror}). For our purposes, we have

\begin{definition}\label{maxunip}
  Given a PVHS of weight $k$ over $(\Delta^*)^r$ whose monodromy is
  unipotent, we say that $0\in \Delta^r$ is a \textsl{maximally
    unipotent boundary point} if
  \begin{enumerate}
  \item $\dim I^{k,k} = 1$, $\dim I^{k-1,k-1} = r$ and $\dim I^{k,k-1}
    = \dim I^{k-2,k}= 0$, where $I^{a,b}$ is the bigrading associated
    to the limiting MHS and,
  \item $\vspan_\C(N_1(I^{k,k}),\ldots,N_r(I^{k,k})) = I^{k-1,k-1}$,
    where $N_j$ are the monodromy logarithms of the variation.
  \end{enumerate}
\end{definition}

Under these conditions, we may identify $\vspan_\C(N_1,\dots,N_r)
\cong \homom(I^{1,1},I^{0,0})$.  Hence, denoting by $\rho\colon
\gp_{-1} \to \homom(I^{1,1},I^{0,0})$ the restriction map, the
subspace $K = \ker(\rho)$ is a canonical complement of
$\vspan_\C(N_1,\ldots,N_r)$ in $\gp_{-1}$.

Note that both, the notion of maximally unipotent boundary point and
the complement $K$ are independent of the choice of basis.

\begin{definition}\label{def:cancoord}
  Let $\VV \to \poly$ be a PVHS having the origin as a maximally
  unipotent boundary point.  A system of local coordinates
  $(q_1,\dots,q_r)$ is called {\sl canonical} if the associated
  holomorphic function $\Gamma_{-1}$ takes values in $K$.
\end{definition}

\begin{prop}\label{prop:special_coordinates}
  Let $\VV \to \poly$ be a PVHS having the origin as a maximally
  unipotent boundary point.  Then there exists, up to scaling, a
  unique system of canonical coordinates.
\end{prop}

\begin{proof}
  Let $s=(s_1,\dots,s_r)$ be an arbitrary system of coordinates around
  $0$.  We can write
  \begin{equation*}
    \rho(\Gamma(s)) \ =\ \sum_{j=1}^r \gamma_j(s)\,N_j,
  \end{equation*}
  where $\gamma_j(s)$ are holomorphic in a neighborhood of $0\in
  \Delta^r$ and $\gamma_j(0)=0$.  The transformation formula
  (\ref{eq:change_w1}) now implies that the coordinate system
  \begin{equation}\label{extension}
    q_j \ :=\ s_j \exp(2\pi i \gamma_j(s))
  \end{equation}
  is canonical.  Moreover, that same formula shows that it is unique
  up to scaling.
\end{proof}

In \cite{ar:Del-local_behavior}, Deligne observed that a variation of
Hodge structure whose limiting MHS is of Hodge-Tate type defines,
together with the monodromy weight filtration, a variation of mixed
Hodge structure.  In this context, the holomorphic functions
$q_{j}(s)$ defined by (\ref{extension}) constitute part of the
extension data of this family of mixed Hodge structures.  He shows,
moreover, that they agree with the special coordinates studied in
\cite{ar:CDGP-pair}, \cite{ar:morrison-guide}, and
\cite{ar:morr-picard-fuchs} for families of Calabi-Yau manifolds in
the vicinity of a maximally unipotent boundary point.  We sketch this
argument for the sake of completeness.

Given a coordinate system $s=(s_1,\dots,s_r)$, let $(W,F_0)$ be the
limiting MHS of weight $k$ and choose $e^0 \in I^{0,0}$.  Let $e^k \in
I^{k,k}$ be such that $Q(e^0,e^k) = (-1)^k$.  Since the origin is a
maximally unipotent boundary point, there exists a basis $e_1^{
  1},\dots,e_r^{ 1}$ of $I^{1 , 1}$ such that $N_j(e_l^{ 1}) =
\delta_{jl} e^0$.  We can define a (multi-valued) holomorphic section
of $\FF^k$ by
\begin{equation*}
  \omega(s)\ :=\ \exp(\sum_{j=1}^r \frac{\log s_j}{2\pi i} N_j)
  \exp(\Gamma(s))\cdot e^k\,.
\end{equation*}
In the geometric setting of a family of varieties $X_s$ the
coefficients
\begin{equation*}
  h_0(s) \ :=\ -Q(e^0, \omega(s)) \quad\hbox{and}\quad h_j(s) \ :=\
  Q(e^1_j, \omega(s))
\end{equation*}
may be interpreted as integrals $\int_\alpha \omega(s)$ over
appropriate cycles $\alpha \in H_k(X_{s_0})$ on the typical fiber.
Clearly, our assumptions imply that $h_0(s) = (-1)^{k+1}$ and
\begin{equation*}
  \begin{split}
    h_j(s)\ &=\ Q(e^1_j, (\sum_l \frac {\log s_l}{2\pi i} N_l +
    \Gamma_{-1}(s))\cdot e^k) \\ &= \ -  \frac {\log s_j}{2\pi i}
    Q(e^0,e^k) - Q(\Gamma_{-1}(s)\cdot e^1_j,e^k)\\ &=\  - \frac {\log
    s_j}{2\pi i} Q(e^0,e^k) - Q(\sum_l\gamma_l(s) N_l e^1_j,e^k)\\ &=\
    (-1)^{k+1} (\frac {\log s_j}{2\pi i} + \gamma_j(s)).
  \end{split}
\end{equation*}
Therefore, $q_j(s) = \exp(2\pi i h_j(s)/h_0(t))$, which agrees with
Morrison's geometric description of the canonical coordinates in
\cite[\S 2]{ar:morr-picard-fuchs}.



\section{Graded Frobenius Algebras and Potentials}
\label{sec:frobenius}

In this section we will abstract the basic properties of the
cup product in the cohomology subalgebra $\oplus_p H^{p,p}(X)$
for a smooth projective variety $X$ and show that this product
structure may be encoded in a single homogeneous polynomial of
degree $3$.

We recall that $(V,*,e_0,\CB)$ is called a Frobenius algebra if
$(V,*)$ is an associative, commutative $\C$-algebra with unit $e_0$,
and $\CB$ is a nondegenerate symmetric bilinear form such that
$\CB(v_1 * v_2, v_3) = \CB(v_1, v_2 * v_3)$.  The algebra is said to
be real if $V$ has a real structure $V_\R$, $e_0\in V_\R$ and both $*$
and $\CB$ are defined over $\R$.  Throughout this paper we will be
interested in {\sl graded}, real Frobenius algebras of weight $k$.  By
this we mean that $V$ has an even grading
\begin{equation*}
  V\ =\ \bigoplus_{p=0}^k V_{2p},
\end{equation*}
defined over $\R$, and such that:
\begin{enumerate}
\item $V_0 \ \cong \ \C$,
\item $(V,*)$ is a graded algebra,
\item $\CB(V_{2p}, V_{2q})  \ =\ 0$ if $p+q \not= k$.
\end{enumerate}

The product structure on a Frobenius algebra $(V,*,e_0,\CB)$ may be
encoded in the trilinear function $\ti{\phi}_0: V\times V \times V
\rightarrow \C$
\begin{equation*}
  \ti{\phi}_0(v_1,v_2,v_3)\ :=\ \CB(e_0,v_1*v_2*v_3),
\end{equation*}
or, after choosing a graded basis $\{e_0,\ldots,e_m\}$ of $V$, in the
associated cubic form:
\begin{equation}\label{eq:phi_0}
  \phi_0(z_0,\ldots,z_m)\ :=\  \frac{1}{6}\ti{\phi}_0(\gamma, \gamma,
  \gamma)\ =\ \frac{1}{6} \CB(e_0, \gamma^3),\quad \gamma \ =\
  \sum_{a=0}^m z_a e_a.
\end{equation}
Let $\{e_a^\CB\ :=\  \sum_b h_{ba}e_b\}$ denote the $\CB$-dual basis of
$\{e_a\}$.

\begin{theorem}\label{thm:prod_pot}
  The cubic form $\phi_0$ defined by~(\ref{eq:phi_0}) is (weighted)
  homogeneous of degree $k$ with respect to the grading defined by
  $\deg z_a = \deg e_a$, and satisfies the algebraic relation:
  \begin{equation}
    \label{eq:prod_pot_2}
    \sum_{d,f} \frac{\del^3 \phi_0}{\del z_a \del z_b \del z_d} h_{df}
    \frac{\del^3 \phi_0}{\del z_f \del z_c \del z_g} \ = \ \sum_{d,f}
    \frac{\del^3 \phi_0}{\del z_b \del z_c \del z_d} h_{df}
    \frac{\del^3 \phi_0}{\del z_a \del z_f \del z_g}.
  \end{equation}
  The product $*$ and the bilinear form $\CB$ may be expressed in
  terms of $\phi_0$ by
  \begin{eqnarray}
    \label{eq:prod_pot_1}
    \CB(e_a, e_b) &=& \frac{\del^3 \phi_0}{\del z_0 \del z_a \del z_b}
    \\
    \label{eq:prod_pot_3}
    e_a * e_b &=& \sum_c \frac{\del^3 \phi_0}{\del z_a \del z_b \del
      z_c} e_c^\CB\,.
  \end{eqnarray}

  Conversely, let $V$ be an evenly graded vector space, $V_0\cong \C$,
  $\{e_0,\ldots, e_m\}$ a graded basis of $V$, and $\phi_0 \in
  \C[z_0,\dots,z_m]$ a cubic form, homogeneous of degree $k$ relative
  to the grading $\deg z_a = \deg e_a$, satisfying
  (\ref{eq:prod_pot_2}).  Then, if the bilinear symmetric form defined
  by (\ref{eq:prod_pot_1}) is non-degenerate, the product
  (\ref{eq:prod_pot_3}) turns $(V,*,e_0,\CB)$ into a graded Frobenius
  algebra of weight $k$.
\end{theorem}

\begin{proof}
  We note, first of all, that the quasi-homogeneity of $\phi_0$
  follows from the assumption that $*$ is a graded product.  Moreover,
  (\ref{eq:prod_pot_1}) is a consequence of the fact that $e_0$ is the
  unit for $*$, while~(\ref{eq:prod_pot_2}) is the associativity
  condition for $*$. On the other hand, (\ref{eq:phi_0}) and the fact
  that $(V,*,e_0,\CB)$ is Frobenius imply that
  \begin{equation*}
    \frac{\del^3 \phi_0}{\del z_a \del z_b \del z_c} \ =\
    \CB(e_0,e_a*e_b*e_c) \ =\ \CB(e_a*e_b,e_c)
  \end{equation*}
  and (\ref{eq:prod_pot_3}) follows.

  The converse is immediate since~(\ref{eq:prod_pot_3}) defines a
  commutative structure whose associativity follows
  from~(\ref{eq:prod_pot_2}); the quasi-homogeneity assumption implies
  that the product is graded; $e_0$ is the unit because
  of~(\ref{eq:prod_pot_1}). The compatibility between $\CB$ and $*$
  comes from $\CB(e_a*e_b, e_c) = \frac{\del^3 \phi_0}{\del z_a \del
    z_b \del z_c} = \CB(e_a, e_b * e_c)$.
\end{proof}

\begin{example}\label{ex:classical_potential_4}
  Let $(V,*,e_0,\CB)$ be a graded real Frobenius algebra of weight
  four.  Setting $\dim V_{2} = r$ and $\dim V_{4} = s$, we have $\dim
  V = 2r + s + 2 :=m + 1$.  We choose a basis $\{T_0,\dots,T_m\}$ of
  $V$ as follows: $T_0 = e_0 \in V_0$ is the multiplicative unit,
  $\{T_1,\dots,T_r\}$ is a real basis of $V_2$, $\{T_{r+1},\dots,
  T_{r+s}\}$ is a basis of $V_4$ such that $B(T_{r+a},T_{r+b}) =
  \delta_{a,b}$.  Finally, $\{T_{r+s+1},\dots, T_{m-1}\}$ and $T_m$
  are chosen as the $\CB$-duals of $\{T_1,\dots,T_r\}$ and $T_0$ in
  $V_6$ and $V_8$, respectively.  We will say that such a basis is
  {\sl adapted} to the graded Frobenius structure.

  It now follows from (\ref{eq:phi_0}), that with respect to such a
  basis, the polynomial $\phi_0$ is given by
  \begin{equation}
    \label{eq:classical_potential_4}
    \phi_0(z) = \frac{1}{2} z_0^2 z_m + z_0 \sum_{j=1}^r z_{j}
    z_{r+s+j}+ \frac{1}{2} z_0 \sum_{a=1}^s z_{r+a}^2 +
    \frac{1}{2} \sum_{a=1}^s z_{r+a} P^a(z_{1},\ldots,z_{r}),
  \end{equation}
  where $P^a(z_{1},\ldots,z_{r}) = \sum_{j,k=1}^r P^a_{jk} z_jz_k$ are
  homogeneous polynomials of degree $2$ determined by:
  \begin{equation}\label{polynomials}
    P^a_{jk} = \CB(T_{r+a}, T_{j} * T_{k}).
  \end{equation}

  Specializing to the case when $V = \oplus_{p=0}^4 H^{p,p}(X)$ for a
  smooth, projective fourfold $X$, endowed with the cup product and
  the intersection form
  \begin{equation*}
    \CB(\alpha,\beta) \ :=\  \int_X \alpha \cup \beta\,,
  \end{equation*}
  we obtain
  \begin{equation*}
    \ti{\phi}_0^{cup}(\omega_1,\omega_2,\omega_3) = \int_X \omega_1
    \cup \omega_2 \cup \omega_3.
  \end{equation*}
\end{example}

In order to complete the analogy with the structure deduced from the
cup product on the cohomology of a smooth projective variety we need
to require that the Lefschetz Theorems be satisfied.  Given $w \in
V_2$, let $L_w \colon V \to V$, denote the multiplication operator
$L_w(v) = w * v$.  The fact that $*$ is associative and commutative
implies that the operators $L_w$, $w\in V_2$, commute, while the
assumption that $*$ is graded implies that $L_w$ is nilpotent and
$L_w^{k+1} = 0$.  Given a basis $w_1,\dots,w_r$ of $V_2$, we can view
$V$ as a module over the polynomial ring $\jgU :=
\C[L_{w_1},\dots,L_{w_r}]$.  Clearly, if $V = \jgU\cdot e_0$, then the
product structure $*$ may be deduced from the action of $\jgU$.  This
conclusion is also true for $k=3$ or $k=4$ without further conditions,
as can be checked explicitly.

If $V$ underlies a real, graded Frobenius algebra of weight $k$, we
can define on $V$ a mixed Hodge structure of Hodge-Tate type by
$I^{p,q} = 0$ if $p\not=q$, and $I^{p,p} = V_{2(k-p)}$.  If we set
\begin{equation}
  \label{eq:Q_from_B}
  Q(v_a, v_b) = (-1)^a \CB(v_a,v_b) \text{ if } v_a \in V_{2a},
\end{equation}
then $Q$ has parity $(-1)^k$ and for $w\in V_2 \cap V_\R$, the
operator $L_w$ is an infinitesimal automorphism of $Q$ and a $(-1,-1)$
morphism of the MHS.  We also observe, that in the geometric case
---such as in Example~\ref{ex:classical_potential_4}--- this
construction agrees with that in Example~\ref{ex:totalcohomology}.

\begin{definition}\label{polarizationforfrobenius}
  An element $w\in V_{2}\cap V_\R$ is said to \textsl{polarize}
  $(V,*,e_0,\CB)$ if $(I^{*,*},Q,L_w)$ is a polarized MHS.  A real,
  graded Frobenius algebra is said to be polarized if it contains a
  polarizing element.  In this case, the cubic form $\phi_0$ is called
  a {\sl classical potential}.
\end{definition}

By Theorem~\ref{th:nilporbit}, a polarizing element $w$ determines a
one dimensional nilpotent orbit $(W(L_w),F)$, where $W_l(L_w) =
\oplus_{2j\geq 2k-l} V_{2j}$ is the weight filtration of $L_w$ and
$F^p = \oplus_{j\leq k-p} V_{2j}$.

Given a polarizing element $w$, the set of polarizing elements is an
open cone in $V_{2} \cap V_\R$. Then, it is possible to choose a basis
$w_1, \ldots, w_r$ of $V_{2}\cap V_\R$ spanning a simplicial cone
${\mathcal C}$ contained in the closure of the polarizing cone and
with $w\in {\mathcal C}$.  Such a choice of a basis of $V_2$ will be
called a {\sl framing} of the polarized Frobenius algebra.

Since the weight filtration is constant over all the elements $L_w$
for $w\in {\mathcal C}$, it follows from Theorem~\ref{th:nilporbit}
that $(W({\mathcal C}),F^*)$ is a nilpotent orbit. Hence, we can
define a polarized VHS on $(\Delta^*)^r$ whose period mapping is given
by
\begin{equation}\label{theta}
  \theta(q_1,\dots,q_r) \ =\ \exp(\sum_{j=1}^r z_j L_{w_j}) \cdot F\
  ;\quad q_j = e^{2\pi i z_j}\,.
\end{equation}
Note that the origin is a maximally unipotent boundary point in the
sense of Definition~\ref{maxunip}.

We conclude this section showing that in the weight-four case,
maximally unipotent, Hodge-Tate, nilpotent orbits yield graded
Frobenius algebras.  In the weight-three case this is done in
\cite[Example~14]{ar:Del-local_behavior}.

\begin{prop}\label{prop:constant_prod_from_nilpotent}
  Let $({N}_1,\ldots, {N}_r;F)$ be a weight-four nilpotent orbit,
  polarized by $Q$, whose limiting MHS is Hodge-Tate.  Suppose that $
  \dim I^{4,4} =1$ and choose a non-zero element ${e_0} \in I^{4,4}
  \cap V_\R$; let $e^*_0 \in I^{0,0} \cap V_\R$ be such that
  $Q(e_0,e_0^*)=1$.  Assume, moreover, that
  $\{N_1(e_0),\dots,N_r(e_0)\}$ are a basis of $ I^{3,3}$.  Let $\CB$
  be obtained from $Q$ as in~(\ref{eq:Q_from_B}).  Then, there exists
  a unique product $*$ on $V$ with unit $e_0$ such that
  \begin{eqnarray}\label{eq:prod_from_NO_1}
    N_j(e_0) * v &:=& N_{j}(v)\,, \text{ for } v\in V,\, j=1,\dots,r
    \\
    \label{eq:prod_from_NO_2}
    v_1 * v_2 &:=& \CB(v_1,v_2)\, e_0^*\,, \text{ for } v_1,v_2\in
    I^{2,2}
  \end{eqnarray}
  Furthermore, $(V,*,e_0,\CB)$ is a graded, polarized, real Frobenius
  algebra.
\end{prop}
\begin{proof}
  It is clear that~(\ref{eq:prod_from_NO_1})
  and~(\ref{eq:prod_from_NO_2}) define a graded product whose unit is
  $e_0$. Commutativity follows immediately from the symmetry of $Q$
  and the commutativity of the operators ${ N}_j$.

  There are two non-trivial cases to check in order to prove the
  associativity of the product. When all three factors lie in
  $I^{3,3}$ this follows, again, from the commutativity of the
  operators $N_j$.  On the other hand, given $v\in I^{2,2}$:
  \begin{equation*}
    \begin{split}
      (N_j(e_0)*N_k(e_0))*v \ &=\ \CB(N_j(e_0)*N_k(e_0),v) \, e_0^*\ =\
      \CB({N}_j({N}_k(e_0)) , v)\, e_0^*\\ &=\ \CB({N}_k({N}_j(e_0)) ,
      v) \, e_0^*\ =\ \CB(e_0, {N}_j({N}_k(v))) \, e_0^*\\ &=\ \CB(e_0,
      N_j(e_0)*(N_k(e_0)*v)) \, e_0^* \\ &=\ N_j(e_0) * (N_k(e_0) *v).
    \end{split}
  \end{equation*}
  Thus, $(V,*, e_0)$ is a graded, commutative, associative algebra
  with unit $e_0$.  It is straightforward to check that $\CB$ is
  compatible with the product.
\end{proof}



\section{Quantum Products}
\label{sec:quantum}

By a quantum product we will mean a suitable deformation of the
(constant) product on a graded, polarized, real Frobenius algebra.
The weight-three case has been extensively studied in the context of
mirror symmetry for Calabi-Yau threefolds (\cite[Chapter
8]{bo:CK-mirror}, \cite{ar:greg-higgs}).  Here we will restrict our
attention to the $k=4$ case.  In order to motivate our definitions, we
recall the construction of the Gromov-Witten potential in the case of
Calabi-Yau fourfolds; we refer to \cite[Ch. 7 and 8]{bo:CK-mirror} for
proofs and details.

As in Example~\ref{ex:classical_potential_4}, let $X$ be a Calabi-Yau
fourfold, and consider the graded, polarized real Frobenius algebra $V
= \oplus_{p=0}^4 H^{p,p}(X)$ , endowed with the cup product and the
intersection form $\CB$.  We choose a basis $\{T_0,\dots,T_m\}$ as in
the example with the added assumption that $\{T_{1}, \ldots, T_{r}\}$
be a $\Z$-basis of $H^{1,1}(X,\Z)$ lying in the closure of the
K\"ahler cone.

Following~\cite[\S~8.2]{bo:CK-mirror}, we define the Gromov-Witten
potential as
\begin{equation}\label{gromovwitten}
  \phi(z) = \phi^{GW}(z) = \sum_n \sum_{\beta\in H_2(X,\Z)}
  \frac{1}{n!} \langle I_{0,n,\beta}\rangle (\gamma^n)q^\beta
\end{equation}
where $\gamma = \sum_{j=0}^m z_j T_j$ and $\langle
I_{0,n,\beta}\rangle$ is the Gromov-Witten invariant
\cite[(7.11)]{bo:CK-mirror}.  The term $q^\beta$ may be interpreted as
a formal power or, given a class $\omega$ in the complexified K\"ahler
cone, as $q^\beta := \exp(2\pi i \int_\beta \omega)$.

The term corresponding to $\beta = 0$ in (\ref{gromovwitten}) yields
the classical potential ${\phi}_0^{cup}(z) = (1/6)\,\int_X \gamma^3$;
moreover, if we set $\delta = \sum_{j=1}^r z_{j} T_{j}$ and $\epsilon
= \gamma - \delta - z_0 T_0$ and apply the Divisor Axiom (see~\cite[\S
8.3.1]{bo:CK-mirror}) we may rewrite (\ref{gromovwitten}) as
\begin{equation*}
  \phi^{GW}(z) =  {\phi}_0^{cup}(z) + \sum_n \sum_{\beta\in
    H_2(X,\Z)-\{0\}} \frac{1}{n!} \langle I_{0,n,\beta}\rangle
  (\epsilon^n) \exp(\int_\beta \delta) q^\beta.
\end{equation*}
Now, the homogeneity properties of the Gromov-Witten potential allow
us to further simplify this expression in case $X$ is a Calabi-Yau
fourfold
\begin{equation*}
  \phi^{GW}(z) = \phi_0^{cup}(z) + \sum_{a=1}^s \sum_{\beta\in
    H_2(X,\Z)-\{0\}} \langle I_{0,1,\beta}\rangle
  (T_{r+a}) z_{r+a} e^{2\pi i \sum_{j=1}^r z_{j} \int_\beta T_{j}}
  q^\beta.
\end{equation*}
Note that the above series depends linearly on $z_{r+1},\dots,z_{r+s}$
while the variables $z_{1},\dots,z_{r}$, appear only in exponential
form. Hence, we can write
\begin{equation*}
  \phi^{GW}(z) = \phi_0^{cup}(z) + \sum_{a=1}^s z_{r+a}
  \phi_h^a(z_{1},\ldots,z_{r}).
\end{equation*}
with $\phi_h^a(z_{1},\ldots,z_{r}) = \Psi^a(e^{2\pi i
  z_{1}},\ldots,e^{2\pi i z_{r}})$. It follows from the Effectivity
Axiom (see~\cite[\S~7.3]{bo:CK-mirror}) that $\Psi^a(0) = 0$.

This construction motivates the following definition of an abstract
potential function for graded, polarized, real Frobenius algebras of
weight four.

\begin{definition}\label{quantumpotential}
  Let $(V,*_0,e_0,\CB)$ be a graded, polarized, real Frobenius algebra
  of weight four and let $\{T_{0},\dots,T_{m}\}$ be an adapted basis as
  in Example~\ref{ex:classical_potential_4}.  Assume, moreover, that
  $T_1,\dots,T_r$ are a framing of $V$.  A \textsl{potential} on
  $(V,*_0,e_0,\CB)$ is a function
  \begin{equation}\label{generalpotential}
    \phi(z)=\phi_0(z) + \phi_\hbar(z) \ ;\quad \phi_\hbar(z) =
    \sum_{a=1}^s z_{r+a}\, \phi_h^a(z_{1},\ldots,z_{r})\,,
  \end{equation}
  where $\phi_0(z)$ is the classical potential associated with
  $(V,*_0,e_0,\CB)$, $\phi_h^a(z_{1},\ldots,z_{r}) =
  \psi^a(q_1,\dots,q_r)$, $q_j =\exp(2\pi i z_{j})$, and $\psi^a(q)$
  are holomorphic functions in a neighborhood of the origin in $\C^r$
  such that $\psi^a(0)=0$.  We will refer to $\phi_\hbar(z)$ as the
  {\sl quantum} part of the potential.
\end{definition}

Given a potential $\phi$ we define a {\sl quantum} product on $V$ by
\begin{equation}\label{quantumproduct}
  T_a * T_b \ =\ \sum_{c=0}^m \frac{\del^3 \phi}{\del z_a \del z_b
    \del z_c} T_c^\CB\,,
\end{equation}
where $\{T_0^\CB,\dots,T_m^\CB\}$ denotes the $\CB$-dual basis.
Clearly, $*$ is commutative and $e_0 = T_0$ is still a unit. The
quantum product is associative if and only if the potential satisfies
the WDVV equations:
\begin{equation*}
  \sum_{a=1}^s \frac{\del^3 \phi}{\del z_{i} \del z_{j} \del z_{r+a}}
  \ \frac{\del^3 \phi}{\del z_{k} \del z_{l} \del z_{r+a}} \
  =\ \sum_{a=1}^s \frac{\del^3 \phi}{\del z_{k} \del z_{j} \del
    z_{r+a}} \ \frac{\del^3 \phi}{\del z_{i} \del z_{l} \del
    z_{r+a}},
\end{equation*}
with $i,j,k,l$ running from $1$ to $r$.  In view of
(\ref{generalpotential}) and (\ref{eq:classical_potential_4}) these
equations are equivalent to (\ref{eq:prod_pot_2}) and
\begin{multline}\label{eq:WDVV2}
  \sum_{a=1}^s (P^a_{ij}\frac{\del^2 \phi_h^a}{\del z_{k} \del z_{l}} +
  P^a_{kl} \frac{\del^2 \phi_h^a}{\del z_{i} \del z_{j}} +
  \frac{\del^2 \phi_h^a}{\del z_{i} \del z_{j}} \frac{\del^2
    \phi_h^a}{\del z_{k} \del z_{l}}) = \\= \sum_{a=1}^s
  (P^a_{il}\frac{\del^2 \phi_h^a}{\del z_{j} \del z_{k}} + P^a_{jk}
  \frac{\del^2 \phi_h^a}{\del z_{i} \del z_{l}} + \frac{\del^2
    \phi_h^a}{\del z_{i} \del z_{l}}\frac{\del^2 \phi_h^a}{\del z_{j}
    \del z_{k}}),
\end{multline}
for all $i,j,k,l=1,\dots,r$, and where $P^a_{ij}$ denotes the
coefficients (\ref{polynomials}).

\begin{remark}  For a classical potential the WDVV equations reduce to the
  algebraic relation (\ref{eq:prod_pot_2}).  The Gromov-Witten
  potential (\ref{gromovwitten}) satisfies the WDVV equations
  (see~\cite[Theorem 8.2.4]{bo:CK-mirror}).
\end{remark}

We can now state and prove the main theorem of this section.
\begin{theorem}\label{main}
  There is a one-to-one correspondence between
  \begin{itemize}
  \item Associative quantum products on a framed Frobenius algebra of
    weight four.
  \item Germs of polarized variations of Hodge structure of weight
    four for which the origin $0\in \C^r$ is a maximally unipotent
    boundary point, and whose limiting mixed Hodge structure is of
    Hodge-Tate type.
  \end{itemize}
  This correspondence, which depends on the choice of an element
  corresponding to the unit, identifies the classical potential with
  the nilpotent orbit of the PVHS while the quantum part of the
  potential is equivalent to the holomorphic function $\Gamma$ defined
  by (\ref{gamma}) relative to a canonical basis.
\end{theorem}

\begin{proof}
  Let $(V,*_0,\CB,e_0)$ be a graded, polarized, real Frobenius algebra
  of weight four.  Let $T_0,\dots,T_m$ be an adapted basis such that
  $T_1,\dots,T_r$ is a framing of $V_2$.  Let $\HH\subset V_2$ be the
  tube domain
  \begin{equation*}
    \HH \ :=\ \{\,\sum_{j=1}^r \,z_j\,T_j\ ;\ {\rm Im}(z_j) > 0\,\}
  \end{equation*}
  We view $\HH$ as the universal covering of $(\Delta^*)^r$ via the
  map $(z_1,\dots,z_r) \mapsto (q_1,\dots,q_r)$, $q_j = \exp(2\pi i
  z_j)$.  Let $\VV$ denote the trivial bundle over $(\Delta^*)^r$ with
  fiber $V$ and $\FF^p$ the trivial subbundle with fiber $\sum_{a \leq
    8-2p} V_{a}$.

  Given a potential $\phi$ and elements $w\in V_2$, $v\in V$, the
  quantum product $w*v$ may be thought of as a $V$-valued function on
  $V$.  Let $w *_s v$ denote its restriction to $\HH$.  It follows
  easily from (\ref{generalpotential}) and
  (\ref{eq:classical_potential_4}) that $w *_s v$ depends only on
  $q_1,\dots,q_r$ and, therefore, it descends to a $V$-valued function
  on $(\Delta^*)^r$, \ie\ a section of $\VV$.  This allows us to
  define a connection $\nabla$ on $\VV$ by
  \begin{equation*}
    \nabla_\pd{}{q_j} v \ := \ \frac 1 {2\pi i q_j} (T_{j} *_s  v)\,,
  \end{equation*}
  where in the left-hand side $v$ represents the constant section
  defined by $v\in V$.  As shown in \cite[Proposition
  8.5.2]{bo:CK-mirror} the WDVV equations for the potential $\phi$
  imply that $\nabla$ is flat.

  We can compute explicitly the connection forms relative to the
  constant frame $\{T_0,\dots,T_m\}$.  We have, for $j,l=1,\ldots,r$
  and $a=1,\ldots,s$:
  \begin{equation}\label{eq:connection}
    \begin{split}
      \nabla_{\pd{}{q_j}} T_0 &= \frac{1}{2\pi i q_j} T_{j}\\
      \nabla_{\pd{}{q_j}} T_{l} &= \frac{1}{2\pi i} \sum_{b=1}^s
      \left( \frac{P^b_{jl}}{q_j} + 2\pi i \pd{}{q_j}(2\pi i q_l
        \pd{\psi^b}{q_l})\right) T_{r+b}\\ \nabla_{\pd{}{q_j}} T_{r+a}
      &= \frac{1}{2\pi i} \sum_{k=1}^r \left(\frac{P^a_{jk}}{q_j} +
        2\pi i \pd{}{q_j}(2\pi i q_k \pd{\psi^a}{q_k})\right)
      T_{r+s+k}\\ \nabla_{\pd{}{q_j}} T_{r+s+l} &= \frac{1}{2\pi i
        q_j} \delta_{jl} T_m\\ \nabla_{\pd{}{q_j}} T_m &= 0,
    \end{split}
  \end{equation}
  where the coefficients $P^a_{jk}$ are defined as in
  (\ref{polynomials}).  Note that these equations imply that the
  bundles $\FF^p$ satisfy the horizontality condition
  (\ref{horizontality}).  Moreover, suppose we define a bilinear form
  $Q$ on $V$ as in (\ref{eq:Q_from_B}) and extend it trivially to a
  form $\QQ$ on $\VV$, then it is straightforward to check that
  \begin{equation*}
    \QQ( \nabla_{\pd{}{q_j}} T_a , T_b) + \QQ( T_a ,
    \nabla_{\pd{}{q_j}}T_b) =0
  \end{equation*}
  for all $j=1,\dots,r$ and all $a,b = 0,\dots,m$.  Hence the form
  $\QQ$ is $\nabla$-flat.

  We may also deduce from (\ref{eq:connection}) that $\nabla$ has a
  simple pole at the origin with residues
  \begin{equation*}
    \res_{q_j=0}(\nabla) \ =\ \frac{1}{2\pi i}\ L^0_{T_j}\ ;\quad
    j=1,\dots,r,
  \end{equation*}
  where $L^0_{T_j}$ denotes multiplication by $T_j$ relative to the
  constant product $*_0$.  It then follows from
  \cite[Th\'eor\`eme~II.1.17]{bo:Del-equations} that, written in terms
  of the (multivalued) $\nabla$-flat basis $T_0^\flat,\dots,T_m^\flat$
  the matrix of the monodromy logarithms $N_j$ coincides with the
  matrix of $2\pi i \res_{q_j=0}(\nabla)$ in the constant basis
  $T_0,\dots,T_m$, \ie\ with $L^0_{T_j}$. Because the operators
  $L^0_{T_j}$ are real, so is the monodromy $\exp(N_j)$ and therefore
  we can define a flat real structure $\VV_\R$ on $\VV$.

  Since $T_1,\dots,T_r$ are a framing of the polarized Frobenius
  algebra $(V,*_0,\CB,e_0)$, it follows from (\ref{theta}) that the map
  \begin{equation*}
    \theta(q_1,\dots,q_r) \ =\ \exp(\sum_{j=1}^r z_j N_{j}) \cdot F\
    ;\quad z_j = e^{2\pi i q_j}\,.
  \end{equation*}
  is the period map of a VHS (a nilpotent orbit) in the bundle
  $(\VV,\VV_\R,\nabla,\QQ)$.  Since the bundles $\FF^p$ are already
  known to satisfy (\ref{horizontality}), we can apply
  Theorem~\ref{th:2.8} to conclude that they define a polarized VHS on
  $(\VV,\VV_\R,\nabla,\QQ)$.

  In order to complete the asymptotic description of the PVHS defined
  by $\FF$ on $\VV$, we need to compute the holomorphic function
  $\Gamma \colon \Delta^r \to \jlg_-$.  Because of
  Theorem~\ref{th:improved_2.8}, it suffices to determine the
  component $\Gamma_{-1}$.  Moreover, it follows from
  Proposition~\ref{prop:special_coordinates} that we may choose
  canonical coordinates $(q_1,\dots,q_r)$ on $\Delta^r$, so that, in
  terms of the basis $T_0,\dots,T_m$, $\Gamma(q)$ has the form:
  \begin{equation}\label{gammacan}
    \Gamma(q) =
    \left(
      \begin{array}{c|c|c|c|c}
        & & & & \\\hline
        0 & & & & \\\hline
        - \transpose{D}(q) & \transpose{C}(q) & & &\\\hline
        *& * & C(q)& &\\\hline
        * & * & D(q) & 0&
      \end{array}
    \right).
  \end{equation}
  Thus, $\Gamma(q)$ is completely determined by the $r\times s$-matrix
  $C(q)$.

  On the other hand, as noted in \S \ref{sec:variations_at_infinity},
  $\Psi(q) = \exp \Gamma(q)\cdot F_0$ is the expression of the Hodge
  bundles $\FF^p$ in terms of the canonical sections
  (\ref{cansections}):
  \begin{equation}\label{cansections2}
    \tilde T(z) \ =\
    \exp (\sum_{j=1}^r\, z_j N_j)\cdot T^{\flat}\,
  \end{equation}
  The matrix $\exp (-\Gamma(q))$, in the basis $T_0,\dots,T_m$, is the
  matrix expressing the canonical sections $\tilde T_0,\dots,\tilde
  T_m$ in terms of the constant frame.  Therefore
  \begin{equation*}
    \tilde T_{r+a}(q) \ =\ T_{r+a} - \sum_{k=1}^r C_{ka}(q) T_{r+s+k}
    - D_a(q)T_m
  \end{equation*}
  and it suffices to compute $\tilde T_{r+a}(q)$ to determine $C$ (and
  $D$).

  It is straightforward to show, using the formulae
  (\ref{eq:connection}), that
  \begin{equation*}
    \begin{split}
      T_m^\flat &= T_m\\ T_{r+s+l}^\flat &= T_{r+s+l} - z_l T_m\\
      T_{r+a}^\flat &=T_{r+a} -\sum_{l=1}^r \pd{(P^a + \phi_h^a)}{z_l}
      T_{r+s+l} + (P^a + \phi_h^a) T_m
    \end{split}
  \end{equation*}
  Hence, to obtain $\tilde T_{r+a}(q)$ it suffices to apply
  (\ref{cansections2}), together with the fact that the matrix of
  $N_j$ in the basis $\{T_p^\flat\}$ coincides with that of
  $*_0$-multiplication by $T_j$ relative to $\{T_p\}$.  Thus,
  $N_j(T_m^\flat) = 0$ and
  \begin{equation*}
    N_j(T_{r+s+l}^\flat)\ =\ \delta_{jl}\, T_m^\flat\,;\quad
    N_j(T_{r+a}^\flat)\ =\ \sum_{k=1}^r P^a_{jk} \,T_{r+s+k}^\flat\,.
  \end{equation*}
  This, together with the fact that $P^a$ is a homogeneous polynomial
  of degree $2$, implies that
  \begin{equation*}
    \begin{split}
      (\sum_{j=1}^r z_j N_j) T^\flat_{r+a} &= \sum_{j,l=1}^r z_j
      \frac{\del^2 P^a}{\del z_j \del z_l} T^\flat_{r+s+l} =
      \sum_{l=1}^r \pd{P^a}{z_l} T^\flat_{r+s+l} \\ &= \sum_{l=1}^r
      \pd{P^a}{z_l} (T_{r+s+l} - z_l T_m) = \sum_{l=1}^r \pd{P^a}{z_l}
      T_{r+s+l} - 2 P^a T_m
    \end{split}
  \end{equation*}
  and
  \begin{equation*}
    \frac{1}{2}(\sum_{j=1}^r z_j N_j)^2 T^\flat_{r+a} =
    \frac{1}{2}(\sum_{j=1}^r z_j N_j) \sum_{l=1}^r \pd{P^a}{z_l}
    T^\flat_{r+s+l} = \frac{1}{2}\sum_{j=1}^r z_j \pd{P^a}{z_j}
    T^\flat_m = P^a T_m.
  \end{equation*}
  Hence
  \begin{equation*}
    \begin{split}
      \ti{T}_{r+a}\ &=\ T^\flat_{r+a} + (\sum_{j=1}^r z_j N_j)
      T^\flat_{r+a} + \frac{1}{2} (\sum_{j=1}^r z_j N_j)^2
      T^\flat_{r+a}\\ &=\ T^\flat_{r+a} + \sum_{l=1}^r \pd{P^a}{z_l}
      T_{r+s+l} - 2 P^a T_m + P^a T_m \\ &=\ T_{r+a} - \sum_{l=1}^r
      \pd{\phi_h^a}{z_l} T_{r+s+l} + \phi_h^a T_m.
    \end{split}
  \end{equation*}
  Thus,
  \begin{equation*}
    C_{ka}\ =\ \pd{\phi_h^a}{z_k}\;\quad \hbox{and}\quad
    D_a\  = \ -\phi_h^a\,.
  \end{equation*}

  Conversely, suppose now that $(\VV,\VV_\R,\FF,\nabla,\QQ)$ is a
  polarized VHS of weight four over $(\Delta^*)^r$, that the origin is
  a maximally unipotent boundary point, and that the limiting MHS is
  of Hodge-Tate type.  Let $\{N_1,\dots,N_r;F\}$ denote the associated
  nilpotent orbit and set $V_{8-2p} := I^{p,p}$.  It follows from
  Proposition~\ref{prop:constant_prod_from_nilpotent} that we can
  define a product $*_0$, and a bilinear form $\CB$ ---as in
  (\ref{eq:Q_from_B})--- turning $(V,*_0,\CB)$ into a polarized, real,
  graded Frobenius algebra with unit $e_0\in V_0$.  This structure is
  determined by the choice of unit and the fact that, relative to an
  adapted basis $\{T_0,\dots,T_m\}$ as in
  Example~\ref{ex:classical_potential_4},
  \begin{equation}\label{coeff}
    N_j (T_{r+a}) \ =\ \sum_{k=1}^r \CB(T_j*_0 T_{r+a}, T_k)\,T_{r+s+k}
    \ =\ \sum_{k=1}^r P^a_{jk}\,T_{r+s+k},
  \end{equation}
  $j=1,\dots,r$, $a=1,\dots,s$, and $P^a_{jk}$ are the coefficients
  (\ref{polynomials}) of the associated classical potential $\phi_0$.

  We have already noted that in canonical coordinates, the holomorphic
  function $\Gamma$ associated with the PVHS takes on the special form
  (\ref{gammacan}).  Moreover, since $\Gamma(q)$ satisfies the
  differential equation (\ref{eq:integcond}), we have from
  (\ref{eq:e7}) that
  \begin{equation*}
    dG_{-2} = \Gamma_{-1} \Theta - \Theta \Gamma_{-1} + \Gamma_{-1}
    d\Gamma_{-1}
  \end{equation*}
  for $\Theta = d(\sum_{j=1}^r z_{j} {N}_j)$. Consequently,
  \begin{equation}\label{int2}
    d(D_a) = - \sum_{k=1}^r C_{ka}\,dz_k.
  \end{equation}

  We now define a potential on $V$ by
  \begin{equation}\label{pot}
    \phi(z)\ :=\  \phi_0(z) - \sum_{a=1}^s z_{r+a} D_a(q)\,.
  \end{equation}
  Since we already know that the classical potential $\phi_0$
  satisfies (\ref{eq:prod_pot_2}), the WDVV equations for $\Phi$
  reduce to the equations (\ref{eq:WDVV2}).  But this is a consequence
  of the integrability condition (\ref{eq:integcond}); indeed, note
  that given (\ref{coeff}), if we let $\Xi = (\xi_{ka})$ be the
  $r\times s$-matrix of one forms
  \begin{equation*}
    \xi_{ka}\ =\ \sum_{j=1}^r (P^a_{jk} + \pd {C_{ka}}{z_j})\,dz_j,
  \end{equation*}
  the equation (\ref{eq:integcond}) reduces to
  \begin{equation}
    \label{eq:WDVV3}
    \Xi \wedge \Xi^t \ =\ 0
  \end{equation}
  which, in view of (\ref{int2}) and (\ref{pot}), is easily seen to be
  equivalent to (\ref{eq:WDVV2}).
\end{proof}

\begin{remark}
  Note that~(\ref{eq:WDVV3}) expresses the WDVV equations in a very
  compact form. Also, note that even though the quantum product is
  defined in terms of third derivatives of $\phi$, one recovers the
  full potential $\phi$ from the PVHS. Since~(\ref{quantumproduct})
  only allows for an ambiguity which is, at most, quadratic in $z$,
  the quantum part $\phi_\hbar$ is uniquely determined.
\end{remark}

\newpage


\providecommand{\bysame}{\leavevmode\hbox to3em{\hrulefill}\thinspace}

\end{document}